\documentclass[conference]{IEEEtran}
\usepackage{cite}

\ifCLASSINFOpdf
   \usepackage[pdftex]{graphicx}
\else
\fi
\usepackage[cmex10]{amsmath}
\usepackage{array}

\begin{document}
%
\title{Some Architectures for Chebyshev Interpolation}

\author{\IEEEauthorblockN{Theja Tulabandhula}
\IEEEauthorblockA{Indian Institute of Technology Kharagpur, India\\ 
Dept. of Electrical Engineering \it{t.theja@iitkgp.ac.in}}
}


%


\maketitle

\begin{abstract}
\small Digital architectures for Chebyshev interpolation are explored and a variation which is word-serial in
nature is proposed. These architectures are contrasted with equispaced system structures. Further, Chebyshev
interpolation scheme is compared to the conventional equispaced interpolation vis-\'{a}-vis reconstruction error
and relative number of samples. It is also shown that the use of a hybrid (or dual) Analog to Digital converter
unit can reduce system power consumption by as much as $1/3^{rd}$ of the original.
\end{abstract}


%
\IEEEpeerreviewmaketitle

\section{Introduction} \label{sec:intro}


Applications like synchronization in software defined radio (SDR) and power constrained sampling in sensor
networks can have solutions garnered from non-uniform sampling research. Often in such pursuits, the hardware
requirements and efficient architecture design are ignored \cite{neagoe90,zhu92}. Signal interpolation is one of
the underlying questions which one tries to solve in such applications. Chebyshev interpolation technique
\cite{fox68} in particular has been a promising non-uniform sampling and interpolation scheme. In general,
sampling on non-uniform grid has many advantages (see \cite{fox68,crc03}). For example, Runge (see \cite{crc03}
pp155-156) demonstrated that interpolation of equispaced signal values is non-optimal for a certain class of
functions. Sampling on the uniform grid, on the other hand though sometimes suboptimal, has been widely used in
clock synchronization, timing correction, sample rate conversion among other applications.

Fox and Parker \cite{fox68} suggested two similar schemes for Chebyshev interpolation. Neagoe et al.
\cite{neagoe90} showed that the coefficient set of one of these interpolation schemes is the output of a DCT
(Discrete Cosine Transform) of the input samples. After these mathematical results Zhu \cite{zhu92}, Wang
\cite{wang96} and Cuypers et al. \cite{cuypers04} have tried presenting digital implementations of the
interpolation scheme. These architectures sometimes don't utilize hardware efficiently or are specific to an
output node set. With an objective of designing a more flexible structure, this paper explores the merits and
demerits of Chebyshev interpolation from an implementation perspective. A systolic array based Chebyshev
interpolation architecture for a window of 8 samples is designed which is word-serial in nature (unlike the
previous suggested structures). A sampling scheme is also proposed involving a SAR (Successive Approximation)
ADC (Analog to Digital Converter) and a flash ADC to make a Flash-SAR hybrid converter block. By suitably
sharing the samples between these ADCs, Chebyshev sampling is performed at $\sim 30-40 \%$ lesser power
consumption levels.

Section \ref{sec:theory} revisits the mathematical basis of Chebyshev interpolation. Digital structures are
explored and a new one is proposed in Section \ref{sec:chebyarch}. Details which make Chebyshev interpolation a
viable alternative to equispaced interpolation are presented in Section \ref{sec:advcheby} and finally Section
\ref{sec:conclusions} summarises the theme and contribution of this paper.




\section{Theory of Chebyshev Interpolation} \label{sec:theory}


Chebyshev polynomials of the first kind $\{T_n(x)\}, x\in [-1,1]$ can be defined recursively as
\begin{equation}
T_{n+1}(x) = 2xT_n(x)-T_{n-1}(x) \label{eqn:cheby_recursive}
\end{equation}
where the first three polynomials are
\begin{eqnarray}
T_0(x) = 1\\
T_1(x) = x\\
T_2(x) = 2x^2-1
\end{eqnarray}

The $k^{th}$ zero of an $n^{th}$ order polynomial ($T_n(x)$) is given as
\begin{equation} \label{eqn:rootchebyshev}
x_k = cos(\frac{2k-1}{2n}\pi),k = 1,2,...,n.
\end{equation}

A polynomial (say $P_N(x)$) can be constructed from $\{T_n(x)\}$, which minimizes the maximum deviation from the
exact underlying signal (\cite{crc03},pp.156). To perform a $N^{th}$ degree polynomial approximation in $[-1,1]$
(this interval can be changed easily), the sample points ($x_k$) should be chosen at the roots of $T_{n+1}(x)$.
This leads to a nonuniform grid which is denser at the edges and sparse towards the center. It can be shown that
the polynomial $P_N(x)$ is a linear combination of $T_0,...,T_{N-1}$ for which the coefficients are the DCT
(Discrete Cosine Transform) of the sample values (sampled at $x_k$) \cite{neagoe90}.

\begin{eqnarray}
P_N(x) = \sum^{N}_{i=0}c_i\bar{T}_i(x) \mbox{ where} \label{eqn:polysum}\\
\bar{T}_0(x) = \frac{1}{\sqrt{2}}T_0(x),\bar{T}_{N>0}(x) = T_{N>0}(x)
\end{eqnarray}

$\{c_j\} \propto$ DCT($[f(x_1) \hdots f(x_{N+1})]^{T}$. We can rewrite this as
\begin{equation}
\left[ \begin{array}{c}
c_0 \\
\vdots \\
c_N \end{array} \right]
=
\mbox{\textbf{C}}
\left[ \begin{array}{c}
f(x_1) \\
\vdots \\
f(x_{N+1}) \end{array} \right] \label{eqn:dctcoeffs}
\end{equation}
with
\begin{equation}
(C)_{j,k} = \mu_j cos\left(\frac{j\pi(2k+1)}{2(N+1)}\right) \label{eqn:cjdctterms0}
\end{equation}
\begin{eqnarray}
\mu_j = \begin{array}{c}
\frac{\sqrt{2}}{N+1},j=0\\
\frac{2(-1)^j}{N+1},j=1,...,N \end{array} \label{eqn:cjdctterms1}
\end{eqnarray}
From Equations \ref{eqn:polysum} and \ref{eqn:dctcoeffs},
\begin{equation}
P_N(x) = [f(x_0)...f(x_N)]\mbox{\textbf{C}}^{T}\mbox{\textbf{T}}
\left[ \begin{array}{c}
x^N\\
\vdots\\
x^0 \end{array} \right] \label{eqn:doublesummation}
\end{equation}

Here \textbf{T} is the matrix of coefficients of powers of $x$ of the Chebyshev polynomials in decreasing order.
In the general Lagrange interpolation case, \textbf{C$^T$T} can be replaced by a matrix \textbf{L} representing
coefficients of Lagrange polynomials.


\section{Architectures for Chebyshev Interpolation}
\label{sec:chebyarch}


\subsection{Prior Art} \label{subsec:priorart}

Two systolic arrays for Chebyshev interpolation by Zhu et al. \cite{zhu92} are
based on transform and time domain descriptions of the interpolation operation respectively. The distinction
between time and transform domain structures is based on which summations in the interpolation formula (Equation
\ref{eqn:doublesummation}) are done first. In the first array, the set of coefficients ${c_i}$ are computed
first and then their product with $\bar{T}_i(x)$ is carried out
. In the second, DCT of the Chebyshev polynomials generates the set of Chebyshev Type Interpolation Functions
(CTIF) $\{\phi_i(x)\}$ which are then used for multiplication with $\{f(x_i)\}$
i.e., Equation \ref{eqn:doublesummation} is rewritten as
\begin{equation}
P_N(x) = \sum_{i=0}^{N}f(x_i)\phi_{i}(x)
\end{equation}
where $\phi_i(x)$ is calculated using Equations \ref{eqn:cjdctterms0} and \ref{eqn:cjdctterms1} as
\begin{equation}
\phi_i(x) = \sum^{N}_{k=0}\mu_k \bar{T}_k(x) cos\left( \frac{k\pi(2i+1)}{2(N+1)}\right)
\end{equation}


A disadvantage with this structure is that the input is assumed to come in parallel. Thus the multiplications
which could have otherwise been scheduled vis-\'{a}-vis time are now being done simultaneously, reducing
hardware utilization efficiency.




A structure similar to the previous ones is proposed by Wang et al \cite{wang96}. It assumes that we output
another set of Chebyshev sampled values (with order $M \neq N$) from the existing ones and the hardware has been
optimized keeping this in view making it unusable for an arbitrary output node set.


Cuypers et al. \cite{cuypers04} also propose two architectures. The first uses a fast DCT block and employs a
fast adder for Chebyshev recursive relations of Equation \ref{eqn:cheby_recursive}. No insight into the
computational load per clock cycle, simultaneous use of adders and other implementation details is given.
The second scheme
 suggests the use of a Farrow structure and a CORDIC unit to perform interpolation
assuming that the input signal is in $\theta$ domain rather than the $x$ domain ($x = cos(\theta)$). Computation
of the Farrow structure coefficients is not explained.

\subsection{Proposed Scheme} \label{subsec:presentart}

\begin{figure*}[htbp]
    \centering
    \includegraphics[scale=0.5]{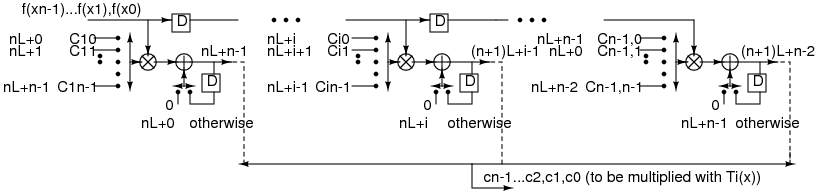}
    \caption{Part of the proposed structure, working on word serial data ${f(x_i)}$ from a data converter to get the coefficients.}
    \label{fig:thejacoeffdigital}
\end{figure*}

The assumption that all the samples are available at the same time is not practical. An implementation which is
word-serial rather than word-parallel in nature would better utilize hardware and not require more buffering of
samples than necessary. Keeping this in view, a design which makes use of the word-serial property of the input
and reduces the overall count of multiply and add units is proposed. Portions of the computation (Equation
\ref{eqn:doublesummation}) are performed as samples arrive one by one. For instance, $\{c_i\}$, the set of
coefficients in the interpolation formula of Equation \ref{eqn:dctcoeffs} is computed using a word serial
systolic array shown in Figure \ref{fig:thejacoeffdigital}. The Chebyshev polynomials are also computed
according to Equation \ref{eqn:cheby_recursive} at an arbitrary node set by rescheduling a pair of multiply and
add units (IIR filtering). Maximum resource usage is guaranteed (i.e., Hardware Usage Efficiency = 100 \%) for
both the computation of the coefficients as well as the subsequent FIR filtering (multiplication and summation).
Though the first part of the total system could be optimized by using any of the fast DCT architectures
available (\cite{parhi05,chiper96}), a generalized systolic array performing matrix vector multiplication has
been used to keep the structure independent of the output node set. The transformation matrix $T$ of the 2D
dependence graph (DG) which was obtained by choosing the desired sequence of inputs and outputs is $T =
[1\hspace{0.2in} 0]$ and the schedule vector ($s$) which allows the reuse of multiply and add units is
$[1\hspace{0.2in}1]^{T}$.

For both, coefficient generation and Chebyshev polynomial evaluation, multiplexers with timing control are used
to achieve correct flow of data. For example, the output of the Chebyshev polynomial evaluation unit has to
switch between IIR filter mode and connect `$0$' and input `$x$' to the output when $T_0(x)$ and $T_1(x)$ are
evaluated in each $N$ cycle period. Note that since the normalized Chebyshev function values are needed, the
recursions are slightly different from Equation \ref{eqn:cheby_recursive}. For example, even though
$\bar{T}_0(x) = \frac{1}{\sqrt{2}}$ and $\bar{T}_{n>0}(x) = T_{n>0}(x)$, the recursion formulae for $\{T_i(x)\}$
will not work for $\{\bar{T}_i(x)\}$. Specifically it will fail at the step where $\bar{T}_2(x)$ is substituted
as $2x\bar{T}_1(x) - \bar{T}_0(x)$ since, $\bar{T}_2(x) = T_2(x) = 2xT_1(x) - T_0(x)$ and $T_0(x)$ is different
from $\bar{T}_0(x)$.

A tabulation of multiplications, additions and computations per cycle required by the word-serial architecture
proposed compared to Zhu's systolic arrays is provided in Table \ref{table:comparison_2}. Solutions by Wang et
al. \cite{wang296} and Cuypers et al. have not been compared because the former optimizes for a specific output
node set and the latter does not discuss the structures at an implementation level.


\begin{table}[htbp]
\addtolength{\tabcolsep}{-2pt} \small
    \centering
    \begin{tabular}{|l|r|r|r|}
    \hline
    Architecture & Zhu                &       Zhu &  Proposed   \\
                 &(time               &(transform &    (1-Dim   \\
                 &  domain)           &   domain) &  systolic)  \\
    \hline
    Buffering required for & samples & $T_i(x)$  & \textbf{none}\\
     (memory)            &\& $T_i(x)$&           &              \\
    \hline
    I/O type    & word              & word      & \textbf{word} \\
                & parallel          & parallel  &\textbf{serial}\\
    \hline
    Computation of coefficients & $\{\Phi_i(x)\}$ & $\{c_i\}$& $\{c_i\}$\\
    Peak operations$/$cycle & $>$8,stored & 8   &  8  \\
    \hline
    Computation of $T_i(x)$ & & & \\
    Peak operations$/$cycle & 0     &   stored  &  8  \\
    \hline
    FIR Filtering ($\sum c_iT_i(x)$)  & & & \\
    Peak operations$/$cycle & 8     &    8      &  8  \\
    \hline
    Latency (cycles)& 16            &   16      &  16 \\
    Hardware Util. Efficiency &100($\%$)&100($\%$)&100($\%$)\\
    \hline
\end{tabular}
\caption{Computations needed in different structures (assuming a set of 8 samples).} \label{table:comparison_2}
\end{table}


\section{Advantages of Chebyshev sampling}
\label{sec:advcheby}


\subsection{Reduction of interpolation error}
\label{subsec:error} 
\begin{figure}[htbp]
    \centering
    \includegraphics[scale=0.34]{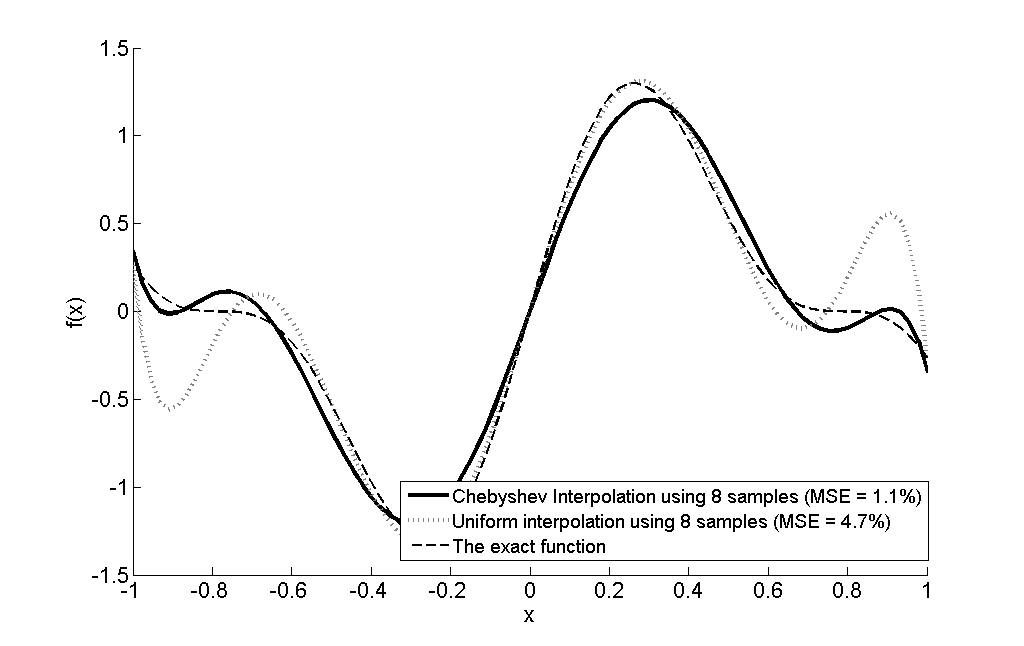}
    \caption{In general, more number of samples for the equispaced case are required in the same interval to reduce the interpolation error.}
    \label{fig:chebyshev_uniform_8_8}
\end{figure}

\begin{figure}[htbp]
    \centering
    \includegraphics[scale=0.34]{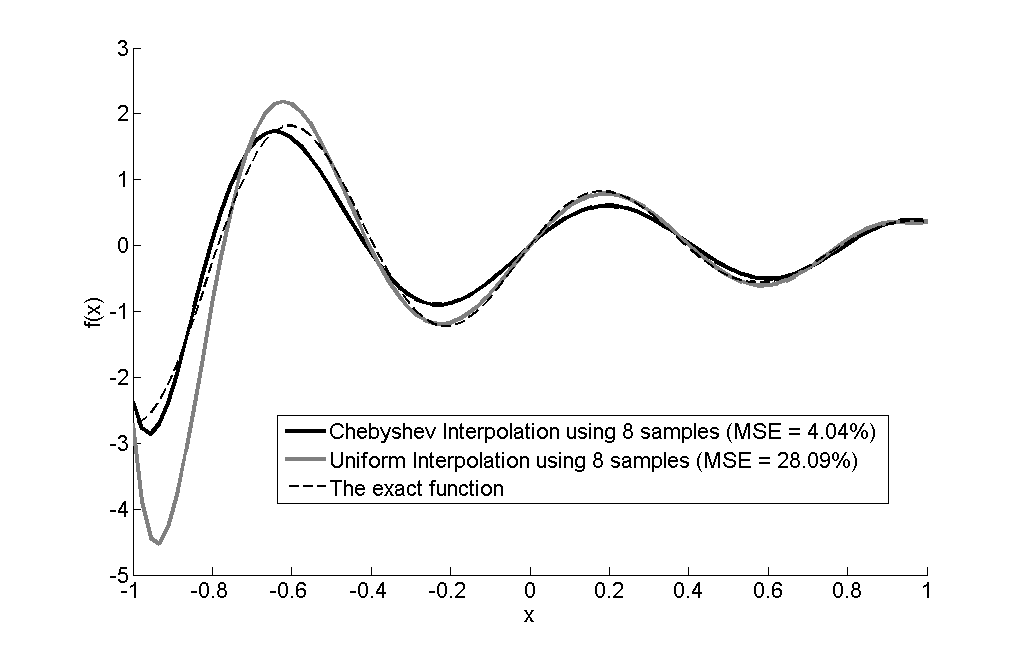}
    \caption{A non bandlimited function $e^{-x}sin(8x)$ over the interval $[-1,1]$ is interpolated using both the schemes.}
    \label{fig:chebyshev_uniform_8_10}
\end{figure}

Two characteristic signals are taken to investigate the effectiveness of the interpolation
scheme. In addition to a bandlimited signal, a non bandlimited signal ($e^{-x}sin(8x)$) is also chosen. When the
number of samples is 8 and the bandlimited function is a sinusoid plus its third harmonic ($sin(4x) +
0.5sin(8x), x \in[-1,1]$), the interpolation error for the equispaced case is $\sim$4 times the Chebyshev case
as shown in Figure \ref{fig:chebyshev_uniform_8_8}. A similar result is obtained for the non bandlimited case as
shown in Figure \ref{fig:chebyshev_uniform_8_10}.

\subsection{Use of Hybrid ADC for power savings} \label{subsec:hybridadc}

Power savings can be achieved during sampling in a Chebyshev based interpolation system through the use of two
(dual) data converters (ADCs) as a hybrid. When the interpolation error limit is fixed for the equispaced and
Chebyshev cases, the number of samples required to do so also becomes fixed. In some cases as seen in Section
\ref{subsec:error},  equispaced system requires more samples than the Chebyshev system. For the Chebyshev
system, a scheme is proposed where the samplings are split between two ADCs, one of which is faster but power
consuming (flash) and the other is slower but power saving (SAR). To do so, a flash and SAR ADC are bundled
together with a timing control unit to make a hybrid unit.
A simple strategy to split the samples between the flash and SAR is based on whether the ratio of the
intersample interval is greater than the SAR sampling and conversion time (i.e., $T_{SAR} <
floor(sin\mbox{kc}/sin\mbox{c})$ where $c=\frac{\pi}{N+1}$ (derivable from Equation \ref{eqn:rootchebyshev}).
Table \ref{table:powersavings} shows the power savings as a function of sharing of samples between the two ADCs
for the two example signals of Section \ref{subsec:error}. Note that, Flash ADC topology is assumed to be
thermometric (not necessarily the case) and power consumption per comparison in either case is taken to be 1
arbitrary unit (au.)



\begin{table}[htbp]
\addtolength{\tabcolsep}{-4pt} \small
    \centering
    \begin{tabular}{|l|r|r|}
    \hline
    ADC unit & 8 bit Flash & 8 bit Flash-SAR hybrid\\
    \hline
    Interpolation & Equispaced & Chebyshev \\
    \hline
    \multicolumn{3}{|l|}{To have interpolation error $<1.1\%$ for $sin(4x)+0.5sin(8x)$ } \\
    \hline
    Num. of points required& 10 & 8 ($a_{flash}+a_{SAR}$)\\
    \hline
    Splitting & --& $a_{flash} = 6$, $a_{SAR} = 2$ \\
    Power (au)& $10*2^8 = 2560$ & $6*2^8 + 2*8 = 1552$ \\
    \hline
    \multicolumn{3}{|l|}{To have interpolation error $<4.1\%$ for $e^{-x}sin(8x)$ } \\
    \hline
    Num. of points required & 11 & 8 ($a_{flash}+a_{SAR}$)\\
    \hline
    Splitting & --& $a_{flash} = 6$, $a_{SAR} = 2$\\
    Power (au)& $11*2^8 = 2816$ & $6*2^8 + 2*8 = 1552$ \\
    \hline
\end{tabular}
\caption{Sampling power savings per window is $\sim39\%$ and $\sim44\%$ resp. for each of the signals.}
\label{table:powersavings}
\end{table}

\subsection{Chebyshev Interpolation and Farrow Structures}
\label{subsec:farrow}

Even though Chebyshev Interpolation using a Farrow Structure is described in \cite{cuypers04}, it assumes that
the input signal is in the $\theta$ domain. It is worthwhile to compare the structures for chebyshev
interpolation in comparison to the Farrow structure which is widely used for equispaced interpolation even
though this kind of polynomial interpolation is ill-conditioned. Farrow structure based interpolation units
recently have been ported to perform some special nonuniform interpolations \cite{babic05} but don't yield to
Chebyshev interpolation because the inter-sample intervals are too diverse in range. The DCT shortcut for the
Chebyshev case has made this scheme comparable to the Farrow shortcut \cite{farrow88} based equispaced case.

\subsection{Design summary and applications} \label{subsec:summaryapplications}


The non-uniformly spaced nodes in the Chebyshev interpolation require block processing which can be a
disadvantage for applications where latency is critical. Further, the sampling times are not only irregular but
they also cause sampling intervals to be non rational ratios of each other. This implies that there will always
be an error in the sampling time even if a very high frequency timing clock is used. An analysis of the optimal
number of sample points to be taken in a Chebyshev window hasn't been done and was fixed to 8. Nevertheless,
this parameter has an effect on the flatness of the system frequency response and on the interpolation error.
From an implementation perspective, latency would increase with an increase in this parameter. In a broader
context like Chebyshev sampling, seeking out optimal node sets contingent to the class of signals at system
input could lead to minimal power consumption and reconstruction errors. But such systems, like Chebyshev
hardware will need extra logic for compatibility with the existing equispaced systems. For Chebyshev sampling to
work well, the average sampling rate should be two times or higher than the $f_{max}$ of the input signal
\cite{neagoe90}, but this is indeed the case in most DSP systems where the equispaced sampling rate is chosen to
be 10-15 times $f_{max}$ as a rule of thumb.

Sampling clock synchronization in DSL modems, timing correction, power efficient sensor networking, sample rate
conversion in Software Defined Radio are some of the topics where Chebyshev interpolation can be used
(fractional delay filters are already being used). Chebyshev interpolation is superior when accuracy is
important. It also fits nicely with signal compression like the DCT compression scheme (where only subset of
coefficients containing most of the energy are retained) \cite{zhu92}.


\section{Conclusions}\label{sec:conclusions}


A digital architecture performing Chebyshev interpolation based on systolic arrays assuming word-serial data
input is implemented in detail and contrasted with other architectures proposed in literature. This structure
has a latency of $2N$ cycles between the input and the interpolated output. Merits of Chebyshev sampling
compared to equispaced sampling are then explored. A scheme for Chebyshev sampling using a Flash-SAR hybrid ADC
unit is also discussed which results in power savings. By optimally distributing the share of samples which will
be sampled by either type of ADC, power consumption of the sampling system is optimized. Once the samples have
been obtained (sequentially), these are fed to the word serial systolic array for interpolation in the digital
domain. Finally, the paper echoes the point that inexpensive digital computation can allow for system specific
optimal node sets (not just Chebyshev and equispaced) leading to arbitrary precision in signal interpolation.

\bibliographystyle{IEEEtran}
\bibliography{chebyREC}

\end{document}